\documentclass[graybox]{svmult}
\usepackage{mathptmx}       
\usepackage{helvet}         
\usepackage{courier}        
\usepackage{type1cm}        
\usepackage{makeidx}         
\usepackage{graphicx}        
\usepackage{multicol}        
\usepackage[bottom]{footmisc}
\usepackage{tikz} 

\usepackage{url}

\usepackage{amsfonts}
\usepackage{amsmath}
\usepackage{amssymb}

\begin{document}

\title{A new Energy Equation Derivation for the Shallow Water Linearized Moment Equations}
\titlerunning{A new Energy Equation Derivation for SWLME}
\author{Julian Koellermeier} 
\institute{Julian Koellermeier \at Department of Mathematics, Computer Science and Statistics, Ghent University \\
Bernoulli Institute, University of Groningen
\at \email{julian.koellermeier@ugent.be}
}  

\maketitle

\abstract{
Shallow Water Moment Equations (SWME) are extensions to the well-known Shallow Water Equations (SWE) for the efficient modeling and numerical simulation of free-surface flows. While the SWE typically assume a depth-averaged vertical velocity profile, the SWME allow for vertical variations of the velocity profile. The SWME therefore assume a polynomial profile and then derive additional evolution equations for the polynomial coefficients via higher order depth integration.
In this work, we perform a new systematic derivation of the energy equation for a specific variant of the SWME, called the Shallow Water Linearized Moment Equations (SWLME). The derivation is based on the standard SWE energy equation derivation and includes the skew-symmetric formulation of the model. The new systematic derivation is beneficial for the extension to other SWME variants and their numerical solution.
}



\section{Shallow Water Equations and energy conservation}
Shallow flows are occurring in a wide range of free-surface flow scenarios, ranging from ocean or river waves to avalanches and mud flows. These shallow flows are characterized by situations in which the surface wave length is much larger than the water depth. For the purpose of this work, we consider a simple 1D setting, i.e., neglecting variations of the geometry or velocity in the lateral (spanwise) direction. Furthermore, we assume a hydrostatic pressure and neglect friction terms. Non-hydrostatic extensions and friction terms can readily be found in \cite{Kowalski,Scholz2024}. Note that the assumption of shallow flows is typically only valid for small bed slopes \cite{Liggett}.

The most commonly used model to simulate such flows is given by the Shallow Water Equations (SWE). They assume a depth-averaged velocity and are derived using depth-averaging of the underlying Navier-Stokes equations. In a 1D setting with known bottom topography term $b$, the system of two PDEs for the water height $h(t,x)$ and the depth-averaged velocity $u_m(t,x)$ reads
\begin{align}
    \partial_t h+\partial_x\left(h u_m\right) = 0, \quad\quad
    \partial_t\left(h u_m\right)+\partial_x\left( h u_m^2 +\frac{1}{2} g h^2\right) =  -gh \partial_x b. \label{eq:SWE_M1}
\end{align}
It is well known that for the SWE \eqref{eq:SWE_M1} 
an energy equation can be derived as \cite{Gassner}
\begin{align}
    \partial_t\left(\frac{h u_m^2}{2}+ g \frac{h^2}{2}+g h b\right) + \partial_x\left(\frac{h u_m^3}{2} +g h u_m (h+b) \right)& =0, \label{eq:SWE_E}
\end{align}
in which the energy $e=\frac{h u_m^2}{2}+ g \frac{h^2}{2}+g h b$ is a conserved quantity and an entropy function with entropy flux $f=\frac{h u_m^3}{2} +g h u_m (h+b)$. This allows tailored numerical methods to preserve energy discretely during numerical simulations \cite{Ersing2025,Gassner,Parisot}.

One way to account for the effects of vertical variations of the horizontal velocity $u(t,x,z)$ over $z$ is to augment the non-linear term $h u_m^2$ in \eqref{eq:SWE_M1} by a so-called Boussinesq coefficient $\beta$ \cite{Liggett}, defined as
\begin{align}
    \beta = \frac{1}{h u_m^2} \int_{b(t,x)}^{h(t,x)} u(t,x,z)^2 \, dz.
    \label{eq:beta}
\end{align}
However, since the precise velocity distribution $u(t,x,z)$ is not known in the SWE \eqref{eq:SWE_M1}, assumptions are needed to model $\beta$, e.g., by assuming a certain shape for $u(t,x,z)$. This ultimately restricts the SWE to only one velocity variable which is the averaged flow velocity $u_m(t,x)$.

\section{Shallow Water Moment Equations}
The standard SWE model \eqref{eq:SWE_M1} is limited by the assumption of depth-averaging, which leads to a single variable $u_m(t,x)$ representing the velocity profile (with the possible inclusion of shape factors or the Boussinesq coefficient $\beta$ \eqref{eq:beta}). This can lead to model inaccuracy in case of velocity profiles that vary significantly in the vertical direction. In fact, this was mentioned in \cite{Liggett} as one of the main causes for errors in open-channel flow calculations. As a way to tackle that, extended shallow water models have been derived. One of these extended models is called Shallow Water Moment Equations (SWME) \cite{Kowalski}.

The SWME uses the vertical velocity expansion
\begin{align}
    u(t,x,z) = u_m(t,x) + \sum_{i=0}^N u_i(t,x) \phi_i\left(\frac{z-b(t,x)}{h(t,x)}\right), \label{eq:expansion}
\end{align}
where $\zeta = \frac{z-b(t,x)}{h(t,x)} \in [0,1]$ is a scaled vertical variable and $u_i(t,x)$ and $\phi_i(\zeta)$ for $i=1,\ldots,N$ are the polynomial coefficients and the orthogonal Legendre polynomials on $[0,1]$, respectively \cite{Kowalski}.

\begin{remark}
    Note that the expansion \eqref{eq:expansion} is the same as a dynamically changing Boussinesq coefficient $\beta$ \eqref{eq:beta}, which becomes a function of the polynomial coefficients $u_i(t,x)$, for $i=1,\ldots,N$ given by
    \begin{align}
        \beta = 1 + \sum_{i=1}^{N} \frac{1}{2i+1} \frac{u_i(t,x)^2}{u_m(t,x)^2} = 1 + \sum_{i=1}^{N} \frac{1}{2i+1} (Mu)_i^2,
    \end{align}
    where $(Mu)_i= \frac{u_i(t,x)}{u_m(t,x)}$ measures the size of the polynomial coefficient $u_i(t,x)$ relative to the depth-averaged velocity $u_m(t,x)$, for $i=1,\ldots,N$. Interestingly, the term $\sum_{i=1}^{N} \frac{1}{2i+1} (Mu)_i^2$ also appears in the computation of Rankine-Hugoniot jump conditions for the below discussed SWLME model in \cite{Koellermeier2022}. This illustrates that the shape of $u(t,x,z)$ has direct impact on properties of the model that can now be modeled explicitly, instead of relying on a single variable $u_m(t,x)$ and an unknown coefficient $\beta$.
\end{remark}

It is clear that $h(t,x), u_m(t,x), u_i(t,x)$ change with time $t$ and space $x$. From now on, we drop the explicit dependence of the variables on $t,x$ to be concise. Furthermore, we assume a given bottom topography $b(t,x) = b(x)$ that is constant in time.

Evolution equations for the variables $h, u_m, u_i$ for $i=1,\ldots,N$ are derived via higher-order depth integration of the Navier-Stokes equations, resulting in $N+2$ equations: one equation for $h$ (continuity equation), one equation for $h u_m$ (depth-averaged momentum balance), and $N$ additional equations for $h u_i$, for $i=1,\ldots, N$ (higher order moment equations), given by
\begin{align}
    \partial_t h+\partial_x\left(h u_m\right) = 0, \quad
    \partial_t\!\left(h u_m\right)+\partial_x\!\left(h u_m^2+ h \sum_{j=1}^N \frac{u_j^2}{2 j+1} +\!\frac{1}{2} g h^2\!\right)\! =  - gh \partial_x b, \label{eq:SWME_M1}\\
    \partial_t\!\left(h u_i\right)+\partial_x\!\left(\!h\!\left(\!2 u_m u_i+\sum_{j, k=1}^N A_{i j k} u_j u_k\!\right)\!\right)\! =  u_m \partial_x\!\left(h u_i\right)\! -\! \sum_{j, k=1}^N \!B_{i j k} u_k \partial_x\left(h u_j\right)\!, \label{eq:SWME_alphai}
\end{align}
with $A_{i j k}=(2 i+1) \int_0^1 \phi_i \phi_j \phi_k d \zeta$ and $B_{i j k}=(2 i+1) \int_0^1 \phi_i^{\prime}\left(\int_0^\zeta \phi_j d \zeta\right) \phi_k d \zeta$, for $i,j,k=1,\ldots,N$, see \cite{Kowalski}.

Different variants of the SWME \eqref{eq:SWME_M1}-\eqref{eq:SWME_alphai} exist by approximating the coefficients $A_{i j k}$ and $B_{i j k}$, see \cite{Koellermeier2022,Koellermeier2020}. As one example, using $A_{i j k}=0, B_{i j k}=0$ leads to the so-called Shallow Water Linearized Moment Equations (SWLME) \cite{Koellermeier2022}. Originally, the SWLME were derived as a simple way to obtain analytical steady states that are extensions of the classical Rankine-Hugoniot jump conditions of the SWE. However, the following energy equation for the SWLME was derived in \cite{Caballero}
\begin{align}
    \partial_t\left(\frac{h u_m^2}{2}+ \frac{ h}{2}\sum\limits_{i=1}^N\frac{u_i^2}{2i+1} + g \frac{h^2}{2}+g h b\right) +\partial_x\left(\frac{h u_m^3}{2} +  \frac{3h u_m}{2} \sum\limits_{i=1}^N \frac{u_i^2}{2i+1} + g h u_m (h+b) \right)&=0.\label{e:SWLME_E1}
\end{align}
This energy equation \eqref{e:SWLME_E1} can be seen as an extension of the standard SWE energy equation \eqref{eq:SWE_E}. We note that the SWLME energy equation \eqref{e:SWLME_E1} takes into account the higher order coefficients $u_i$ for $i=1,\ldots,N$ in the energy.

The derivation of the energy equation \eqref{e:SWLME_E1} in \cite{Caballero} was only a side remark and left out several intermediate steps, that could inform the generalization to more models of similar type. In this paper, we therefore give a new derivation based on an extension of the step-by-step derivation for the standard SWE model given in \cite{Gassner}.

\section{Energy Equation Derivation for the SWLME}
\label{sec:SWLME}
The frictionless 1D SWLME system with bottom bathymetry as derived in \cite{Koellermeier2022}, consists of $N+2$ equations for $h$, $hu_m$ and $h u_i$ for $i=1, \ldots, N$ and is written as
\begin{align}
    (\textrm{C}):&&\partial_t h+\partial_x\left(h u_m\right) &= 0, \label{eq:SWLME_C}\\
    (\textrm{M}):&&\partial_t \left(h u_m\right)+\partial_x \left(h u_m^2+ h \sum_{j=1}^N \frac{u_j^2}{2 j+1} + \frac{1}{2} g h^2 \right)  &=  - gh \partial_x b, \label{eq:SWLME_M1}\\
    (u_i):&&\partial_t \left(h u_i\right) + \partial_x \left( 2 h u_m u_i \right) &=  u_m \partial_x \left(h u_i\right), \label{eq:SWLME_alphai}
\end{align}
with (C) indicating the continuity equation, (M) the momentum balance, and ($u_i$) the moment coefficient $i$ balance, for $i=1,\ldots, N$.

Note again that the SWLME is a special case of the SWME \eqref{eq:SWME_M1}-\eqref{eq:SWME_alphai} as it uses the coefficients
\begin{equation}\label{eq:Aijk_Bijk_SWLME}
    A_{i j k}=0, \quad B_{i j k}=0.
\end{equation}

We will now derive an energy equation for the SWLME \eqref{eq:SWLME_C}-\eqref{eq:SWLME_alphai} step-by-step.\\

Firstly, the potential energy $(\textrm{P})$ is obtained in exactly the same way as for the SWE \cite{Gassner} by computing $g(h+b) \cdot (\textrm{C}) = (\textrm{P})$, which results in
\begin{align}
    (\textrm{P}):&&  \partial_t\left(g \frac{h^2}{2}+g h b\right)+g(h+b) \partial_x(h u_m)& =0, \label{eq:SWLME_P}
\end{align}
where $g \frac{h^2}{2}+g h b$ denotes the potential energy.

Secondly, we aim to derive an equation for the (averaged) kinetic energy.
We notice that the momentum balance (M) can be written as
\begin{align}
    (\textrm{M'}):&&\partial_t\left(h u_m\right)+\partial_x\left( h u_m^2 + h \sum_{j=1}^N \frac{u_j^2}{2 j+1}\right) + gh \partial_x (h+b) &=0. \label{eq:SWLME_M2}
\end{align}

We can then compute $(\textrm{M'})- u_m\cdot (\textrm{C})$ to obtain the new equation (\textrm{A}) as
\begin{align}
    &&\hspace{-1cm} \underbrace{\partial_t (hu_m) - u_m \partial_t h} + \underbrace{\partial_x(h u_m^2) - u_m \partial_x (hu_m)} + gh \partial_x (h+b) + \partial_x  \sum_{j=1}^N \frac{hu_j^2}{2 j+1}&= 0,\\
    (\textrm{A}):&& h \partial_t u_m \quad\quad+ \quad\quad\quad hu_m \partial_x u_m \quad\quad+gh \partial_x (h+b) + \partial_x  \sum_{j=1}^N \frac{hu_j^2}{2 j+1}&= 0. \label{eq:SWLME_A}
\end{align}

Next, we average (A) and (M') to obtain $(\textrm{S}) = \frac{1}{2} (\textrm{A}) + \frac{1}{2} (\textrm{M'})$
\begin{align}
    (\textrm{S}):&&\frac{1}{2}\left(\partial_t(h u_m)+h \partial_t u_m\right) + \frac{1}{2}\left(\partial_x\left(hu_m^2\right)+ h u_m \partial_x u_m\right)+g h \partial_x(h+b) + \partial_x  \sum_{j=1}^N \frac{hu_j^2}{2 j+1}&=0. \label{eq:SWLME_S}
\end{align}
Note that $(\textrm{S})$ is skew-symmetric in the $u_m$-derivatives, see \cite{Gassner}.

We can then compute $u_m \cdot (\textrm{S})$ to obtain an equation for the kinetic energy (K)
\begin{align}
    &&\hspace{-1cm} \frac{1}{2}(\underbrace{ u_m \partial_t(h u_m)+h u_m \partial_t u_m})+\frac{1}{2}(\underbrace{u_m \partial_x\left(h u_m^2\right)+h u_m^2 \partial_x u_m})+g h u_m \partial_x(h+b)+ u_m \partial_x  \sum_{j=1}^N \frac{hu_j^2}{2 j+1}&=0\\
    (\textrm{K}):&& \partial_t\left(\frac{h u_m^2}{2}\right)\quad\quad +  \quad\quad\quad \partial_x\left(\frac{h u_m^3}{2}\right) \quad\quad\quad+g h u_m \partial_x(h+b)+ u_m \partial_x  \sum_{j=1}^N \frac{hu_j^2}{2 j+1}& =0, \label{eq:SWLME_K}
\end{align}
where $\frac{hu_m^2}{2}$ denotes the kinetic energy.

Thirdly, we compute an equation from the moment equations ($u_i$) by following the same recipe as above, only applied to ($u_i$) instead of (M). We start by computing the alternative moment equation ($u_i'$) using elementary changes as
\begin{align}
    (u_i'):&& \partial_t\!\left(h u_i\right) + \partial_x\!\left(\! h u_m u_i \right) + h u_i \partial_x u_m &=  0. \label{eq:SWLME_alphai2}
\end{align}
Then we compute $(\textrm{A$u_i$}) = (\textrm{$u_i'$})- u_i\cdot (\textrm{C})$ as
\begin{align}
    &&\underbrace{\partial_t (h u_i) - u_i \partial_t h} + \underbrace{\partial_x(h u_m u_i) - u_i \partial_x(h u_m)} + hu_i \partial_x u_m&= 0,\\
    (\textrm{A$u_i$}):&& h \partial_t u_i \quad\quad+ \quad\quad\quad hu_m \partial_x  u_i \quad\quad\quad + hu_i \partial_x u_m &= 0, 
\end{align}

Next, we again average $(\textrm{A}u_i)$ and $(u_i')$ to obtain $(\textrm{S}u_i) = \frac{1}{2} (\textrm{A}u_i) + \frac{1}{2} (u_i')$
\begin{align}
    (\textrm{S}u_i):&&\frac{1}{2}\left(\partial_t(h u_i)+h \partial_t u_i\right) + \frac{1}{2}\left( \partial_x\left(hu_m u_i\right)+ h u_m \partial_x u_i\right) + h u_i \partial_x u_m &=0. \label{eq:SWLME_Salphai}
\end{align}
Note that similar to before $(\textrm{S}u_i)$ is skew-symmetric in the $u_i$-derivatives, see \cite{Gassner}.

We then compute $u_i \cdot (\textrm{S}u_i)$ to get
\begin{align}
    &&\frac{1}{2}\left(\underbrace{u_i\partial_t(h u_i)+h u_i \partial_t u_i}\right) + \frac{1}{2}\left(\underbrace{ u_i \partial_x\left(hu_m u_i\right)+ h u_m u_i \partial_x u_i}\right) + h u_i^2 \partial_x u_m &=0,\\
    && \partial_t\left(\frac{h u_i^2}{2}\right)\quad\quad +  \quad\quad\quad\quad \partial_x\left(\frac{h u_m u_i^2}{2}\right) \quad\quad\quad + h u_i^2 \partial_x u_m & =0, 
\end{align}

or equivalently
\begin{align}
    (\textrm{K}u_i):&& \partial_t  \left( \frac{h}{2}\frac{u_i^2}{2i+1}\right) +  \partial_x \left(  \frac{h u_m}{2}\frac{u_i^2}{2i+1}\right) + \frac{h u_i^2}{2i+1} \partial_x u_m  &= 0, \label{eq:SWLME_Aalphai}
\end{align}
in which $ \frac{h}{2}\frac{u_i^2}{2i+1}$ can be seen as a "partial kinetic energy" including the effect of the moment coefficients $u_i$.

We can then compute an equation for the total kinetic energy $(\textrm{K}u)$ using the formula $(\textrm{K}u) = (\textrm{K}) + \sum\limits_{i=1}^N (\textrm{K}u_i)$
\begin{align}
    &&\partial_t\left(\frac{h u_m^2}{2}+ \frac{ h}{2}\sum\limits_{i=1}^N\frac{u_i^2}{2i+1}\right) + \partial_x \frac{h u_m^3}{2} +g h u_m \partial_x(h+b) \quad\quad\quad\quad\quad\quad&  \nonumber \\
    &&+ \underbrace{u_m \partial_x  \sum_{j=1}^N \frac{hu_j^2}{2 j+1} + \sum\limits_{i=1}^N\frac{u_i^2 h}{2i+1} \partial_x u_m +  \partial_x \left(  \frac{h u_m}{2} \sum\limits_{i=1}^N \frac{u_i^2}{2i+1}\right)}& =0, \\
    &&\partial_x \left(  \frac{3 h u_m}{2} \sum\limits_{i=1}^N \frac{u_i^2}{2i+1}\right)\quad\quad\quad\quad\quad\quad\quad\quad& \nonumber
\end{align}
which means that we have for the total kinetic energy $(\textrm{K$u$})$:
\begin{align}
    \partial_t\left(\frac{h u_m^2}{2}+ \frac{ h}{2}\sum\limits_{i=1}^N\frac{u_i^2}{2i+1}\right) + \partial_x \frac{h u_m^3}{2} +g h u_m \partial_x(h+b) + \partial_x \left(  \frac{3h u_m}{2} \sum\limits_{i=1}^N \frac{u_i^2}{2i+1}\right)& =0, \label{eq:SWLME_Ktotal}
\end{align}
where $\frac{hu_m^2}{2} + \frac{ h}{2}\sum\limits_{i=1}^N\frac{u_i^2}{2i+1}$ denotes the total kinetic energy.

By adding the equations for the total kinetic energy (K$u$) and the potential energy (P), we obtain an equation for the total energy (E):
\begin{align}
    \! \partial_t\left(\frac{h u_m^2}{2}+ \frac{ h}{2}\sum\limits_{i=1}^N\frac{u_i^2}{2i+1} + g \frac{h^2}{2}+g h b\right) +\partial_x\left(\frac{h u_m^3}{2} +  \frac{3 h u_m}{2} \sum\limits_{i=1}^N \frac{u_i^2}{2i+1} + g h u_m (h+b) \right)&=0, \label{eq:SWLME_E}
\end{align}
where we denote the total energy $e=\frac{h u_m^2}{2}+ \frac{ h}{2}\sum\limits_{i=1}^N\frac{u_i^2}{2i+1}+ g \frac{h^2}{2}+g h b$. It can directly be seen that the total energy $e$ is a conserved quantity and an entropy function with entropy flux $f=\frac{h u_m^3}{2} +  \frac{3 h u_m}{2} \sum\limits_{i=1}^N \frac{u_i^2}{2i+1} + g h u_m (h+b)$.

The energy $e$ can be written in convective variables $(h,q,r_i)=(h, hu_m, hu_i)$ as
\begin{align}
    e&=\frac{q^2}{2h}+ \frac{ 1}{2h}\sum\limits_{i=1}^N\frac{r_i^2}{2i+1}+ g \frac{h^2}{2}+g h b.
\end{align}
This means that the associated entropy variables can readily be computed as
\begin{align}
    q_1&=\partial_h e = -\frac{u_m^2}{2}- \frac{1}{2}\sum\limits_{i=1}^N\frac{u_i^2}{2i+1}+g(h+b), \label{eq:SWLME_q1}\\
    q_2&=\partial_{q} e = u_m, \label{eq:SWLME_q2}\\
    q_{u_i}&=\partial_{r_i} e = \frac{u_i}{2i+1}. \label{eq:SWLME_qalphai}
\end{align}

We note that the total energy equation can therefore equivalently be obtained by computing $ (E) = q_1 \cdot (\textrm{C}) + q_2 \cdot (\textrm{M}) + \sum\limits_{i=1}^N q_{u_i} \cdot (\textrm{$u_i$})$.

\section{Conclusion}
This paper presented a new systematic energy equation derivation for the SWLME model of free-surface flows.

The SWLME energy equation \eqref{eq:SWLME_E} is the same as \eqref{e:SWLME_E1} derived in \cite{Caballero}. However, the derivation presented in this paper is based on an extension of the SWE energy equation derivation in \cite{Gassner}. It reveals as a byproduct a skew-symmetric form of both the momentum equation \eqref{eq:SWLME_S} and the moment equations \eqref{eq:SWLME_Salphai}.

While the derived energy equation \eqref{eq:SWLME_E} is specific to the SWLME system \eqref{eq:SWLME_C}-\eqref{eq:SWLME_alphai}, the derivation of the energy follows a clear set of steps. These steps are an extension of the energy equation derivation for the standard SWE in \cite{Gassner}. We therefore believe that the derivation in this paper is useful for the derivation of energy equations for more generalized systems like the SWME \cite{Kowalski} or dispersive models \cite{Scholz2024}. In addition, the skew-symmetric formulation of the problem could be used to derive entropy conservative numerical schemes similar to \cite{Bohm2020,Ersing2024,Ersing2025,Gassner}.

\begin{acknowledgement}
This publication is part of the project \textit{HiWAVE} with file number VI.Vidi.233.066 of the \textit{ENW Vidi} research programme, funded by the \textit{Dutch Research Council (NWO)}.
\end{acknowledgement}


\end{document}